\documentclass[11pt]{article}
\usepackage{amsfonts,amsmath,amssymb,amsthm, amscd}
\usepackage[dvips]{epsfig,graphics}

\numberwithin{equation}{section}
\newtheorem{theorem}{Theorem}[section]
\newtheorem{prop}[theorem]{Proposition}
\newtheorem{lemma}[theorem]{Lemma}
\newtheorem{cor}[theorem]{Corollary}
\newtheorem{conj}[theorem]{Conjecture}

\theoremstyle{definition}
\newtheorem{definition}[theorem]{Definition}
\newtheorem{example}[theorem]{Example}
\newtheorem{remark}[theorem]{Remark}
\newtheorem{question}[theorem]{Question}

\def\<{{\langle}}
\def\>{{\rangle}}

\def\a{{\alpha}}

\def\e{{\epsilon}}
\def\g{{\gamma}}
\def\S{\mathbb S}

\def\a{\alpha}
\def\s{\sigma}

\def\e{\epsilon}

\def\bs{\bigskip}
\def\Z{\mathbb Z}
\def\C{\mathbb C}

\def\N{[1,N]}

\def\D{{\Delta}}

\begin{document}

\title{Twisting Alexander Invariants with Periodic Representations}

\author{Daniel S. Silver \and Susan G. Williams\thanks{Both authors partially supported by NSF grant
DMS-0706798.} \\ {\em
{\small Department of Mathematics and Statistics, University of South Alabama}}}

\maketitle %{\setlength{\linewidth}{2in}

\begin{abstract} Twisted Alexander invariants were previously defined for any knot and linear representation of its group $\pi$. The invariants are generalized by using any periodic representation of the commutator subgroup $\pi'$. Properties of the new twisted invariants are given. Under suitable hypotheses, reciprocality and bounds on the moduli of zeros are obtained.
A topological interpretation of the Mahler measure of the invariants is presented.    


\end{abstract} 

%\noindent {\it Keywords:} Knot, twisted Alexander polynomial, representation shift, Mahler measure\begin{footnote}{Mathematics Subject Classification 2010:  
%Primary 57M25; secondary 37B10, 37B40.}\end{footnote}

\section{Introduction} The Alexander polynomial was a product of J.W. Alexander's efforts in the 1920s to compute orders of torsion subgroups  of $H_1 M_r$ for various $r$-fold cyclic covers $M_r$ of a knot $k$.  Twisted Alexander polynomials, introduced in 1990 by X.-S. Lin \cite{lin01}, are more sensitive, employing information from linear representations of the knot group $ \pi_1(\S^3 \setminus k)$. Later work by M. Wada \cite{wada94} showed that the twisted Alexander invariants, like their untwisted ancestors, can be defined for finitely presented groups having infinite abelianization.

Any knot group $\pi$ is finitely generated, and hence it admits only finitely many representations into a given symmetric group $S_N$. However, the commutator subgroup $\pi'= [\pi, \pi]$ often admits uncountably many. In \cite{swTAMS} the authors showed that the set of such representations has the structure of a shift of finite type, a dynamical system with attractive properties. Periodic elements of these {\it representation shifts} carry information about the various branched cyclic covers $M_r$. 

Here we show that for any periodic element $\rho$ of ${\rm Hom}(\pi', S_N)$, a twisted Alexander polynomial is defined. Much of our motivation comes from comments on page 637 of \cite{kl99}. 

In most cases, $\rho$ does not extend over the knot group $\pi$, and the resulting invariant is new. When it does extend, the polynomial we obtain is essentially Lin's polynomial. Regardless, there are computational advantages to working with $\pi'$ rather than $\pi$. For example, metabelian representations of $\pi$ restrict to more manageable, abelian representations of $\pi'$; we exploit this in the last section to obtain bounds on the zeros of twisted Alexander polynomials when finite metabelian representations are used.

We formulate our results for general finitely presented groups whenever possible. As is the case for Wada's invariant, more can be proved for the fundamental groups of 3-manifolds. 

The authors are grateful to Stefan Friedl for his suggestions and helpful comments. 

\section{Twisted Alexander invariants} 

\begin{definition} An {\it augmented group system} is a triple ${\cal G}=(G, \e, x)$ consisting of a finitely presented group G, epimorphism  $\e: G \to \<t \mid \> \cong \Z$, and an element $x \in G$ such that $\e(x) = t$. \end{definition}

Augmented group systems were introduced in \cite{silver93}. They are objects of a category in which a morphism from ${\cal G}_1 = (G_1, \e_1, x_1)$ to ${\cal G}_2 = (G_2, \e_2, x_2)$ is a homomorphism $h: G_1 \to G_2$ such that $\e_2 \circ h = \e_1$ and $h(x_1) = x_2$. 

\begin{definition} The augmented group system of a knot $k$ is the triple ${\cal G} =(G, \e, x)$ where $G =\pi_1(\S^3 \setminus k)$, $\e: G \to \Z$ is abelianization, and $x$ is a meridian. \end{definition}

Up to isomorphism, ${\cal G}$ is independent of the choice of meridian $x$. 

A pair $(G, \e)$ is called an ${\it augmented\ group}$. Such pairs are also objects of a category, with morphisms defined analogously. 

Assume that $R$ is a Noetherian unique factorization domain. (Here we are primarily concerned with $R=\Z$.) Let $\rho: G \to {\rm GL}_N R$ be a linear representation. Define 
$$\e \otimes \rho: \Z[G] \to M_N(R[t^{\pm 1}])$$
by mapping $g\in G $ to $\e(g)\rho(g)$ and extending linearly over the group ring $\Z[G]$. 

Let ${\cal G}=(G, \e, x)$ be any augmented group system, and let $ \< x, y, \ldots, z \mid r, \ldots, s\>$ be a presentation of $G$. We assume without loss of generality that the number of generators exceeds the number of relators. 

Fox derivatives  
of words in  the free group $F = \<x, y, \ldots, z \mid\>$ are elements of the group ring $\Z[F]$.   The natural projection $F  \to G$ extends to a ring homomorphism 
$\phi: \Z[F] \to \Z[G]$ that we use to define the 
the {\it Jacobian matrix}:
\begin{equation}\label{jacobian} J=  \begin{pmatrix} \phi({{\partial r}\over {\partial y}}) & \dots &  \phi({{\partial r}\over {\partial z}}) \\ \vdots &&\vdots\\ \phi({{\partial s}\over {\partial y}}) & \dots &  \phi({{\partial s}\over {\partial z}})\end{pmatrix}.\end{equation}
The {\it $\rho$-twisted Jacobian matrix} $J_{\rho}$ is obtained by evaluating each entry of $J$ by $\e\otimes \rho$, and then by removing the inner parentheses of the resulting block matrix. The entries of $J_\rho$ are in $R[t^{\pm 1}]$.

\begin{definition} The {\it Alexander-Lin polynomial} $D_\rho(t)$ is the order $\D_0(J_\rho)$ (i.e., the greatest common divisor of the maximal minors of $J_\rho$) defined up to multiplication by units in $R[t^{\pm 1}]$. \end{definition}

Alexander-Lin polynomials, motivated by Lin's original definition of twisted Alexander polynomials, were introduced in \cite{swALP}. The polynomial $D_\rho(t)$ depends only on the augmented group system ${\cal G}= (G, \e, x)$ and representation $\rho$. The proof of the following is left to the reader. 

\begin{lemma} If there exists $S\in {\rm GL}_NR$ such that $\rho'(g) = S \rho(g) S^{-1}$ for all $g \in G$, then $D_{\rho'}(t) = D_\rho(t)$. \end{lemma} 

\begin{cor} $D_\rho(t)$ depends on $x$ only up to conjugacy in $G$. \end{cor} 

\section{Representation shifts} As in the previous section, assume that ${\cal G} = (G, \e, x)$ is an augmented group system. Denote the kernel of $\e$ by $K$.  Given a presentation $$G= \<x, y, \ldots, z \mid r, \ldots, s\>,$$ we can easily convert it into a presentation of the form 
\begin{equation} \label{presentG} G= \<x, a, b, \ldots c \mid r, \ldots, s\>,\end{equation}
where $a, b, \ldots, c \in K.$ (This can be done by introducing new generators $a, b, \dots, c$ with defining relations $a = yx^{-\e(y)}, \ldots, c = zx^{-\e(z)}$, and then using the new relations to eliminate $y, \ldots, z$ from the presentation.) 

The Reidemeister-Schreier method (see \cite{ls}, for example) enables us to write a presentation for $K$. 

\begin{equation}\label{presentK} K= \< a_\nu, b_\nu, \ldots, c_\nu  \mid r_\nu, \ldots, s_\nu \>.\end{equation}
Here $\nu$ ranges over the integers. The symbol $a_\nu$ is an abbreviation for $x^\nu a x^{-\nu}$ with similar meanings for  $b_\nu, \ldots, c_\nu$. Likewise, $r_\nu, \ldots, s_\nu$ represent $x^\nu r x^{-\nu}, \ldots, x^{\nu} s x^{-\nu}$, rewritten in terms of $a_\nu, b_\nu, \ldots, c_\nu$. 

We regard the symmetric group $S_N$ on symbols $\N=\{1, \ldots, N\}$ as a group of permutation matrices, elements of ${\rm GL}_N\Z$. 
A representation $\rho: K \to S_N$ is an assignment of free group generators $a_\nu, \ldots, c_\nu$ such that the words $r_\nu, \ldots, s_\nu$ are mapped to the identity. \bs

{\sl We denote by $\mu$ the automorphism of $K$ given by $u \mapsto xu x^{-1}$, for all $u \in K$.} 

\begin{definition} If $\rho: K \to S_N$ is a representation, then $\s\rho: K \to S_N$ is the representation $\rho \circ \mu$; i.e., 
$(\s \rho)(u) = \rho(x u x^{-1})$, for all $u \in K$. \end{definition} 

The authors showed in \cite{swTAMS} that ${\rm Hom}(K, S_N)$ has the structure of a compact metric space $\Phi_{S_N}({\cal G})$ and, with this structure, $\s$ is a self-homeomorphism. Elements $\rho$ correspond to the bi-infinite paths of a finite directed graph. Of special interest are the cycles. 

\begin{definition} $\rho \in \Phi_{S_N}({\cal G})$ is {\it periodic} (with {\it period} $r$) if $\s^r \rho = \rho$. \end{definition} 

\section{Crowell group of $K$}  Given a group and permutation representation, R.H. Crowell introduced the {\it derived group} in \cite{crowellderived}. His motivation was to unify algebraically the notions of knot group, Alexander matrix and covering space. While elements of Crowell's derived groups can be interpreted in several ways, the most natural from the topologist's perspective is as lifted paths in covering spaces. 

The derived group $G_\rho$ of a knot group was used in \cite{swCDG} to recover twisted Alexander invariants. The derived group of the commutator subgroup $K_\rho$, which we consider below, is simpler.

Let ${\cal G} = (G, \e, x)$ be an augmented group system. Assume that $\rho \in \Phi_{S_N}({\cal G})$. 
The group $K$ acts on $\{1, \ldots, N\}$ on the right via $\rho$. The {\it Crowell group $K_\rho$} has generator set $\N \times K$, which we will write as  $\{{}^iu \mid i \in \N, u \in K\}$, and relations 
${}^i(uv) = ({}^iu)({}^{i\rho(u)}\!v)$, for all $u, v \in K, 1\le i \le N$. (When $\rho$ is understood, we abbreviate $i \rho(u)$ by $i u$.) It is easily seen that if $K$ has presentation (\ref{presentK}), then
\begin{equation} \label{presentKrho} K_\rho = \< {}^i a_\nu, {}^ib_\nu, \ldots, {}^ic_\nu \mid {}^ir_\nu, \ldots, {}^is_\nu \>.\end{equation} Here $i$ ranges over $\N$, while $\nu$ ranges over $\Z$. Each relation can be expressed in terms of generators using ${}^i(uv) = ({}^iu)({}^{iu}\!v)$. 

The following proposition follows immediately from the definition of $K_\rho$.

\begin{prop} \label{quotient} For any $\rho \in \Phi_{S_N}({\cal G})$, the mapping ${}^iu \mapsto u, \ u \in K$ induces a projection of $p: K_\rho \to K$.

\end{prop}

The Crowell group $K_\rho$ is an invariant of the orbit class of $\rho$, as the next proposition shows. 

\begin{prop} \label{iso} For any $\rho \in \Phi_{S_N}({\cal G})$, the Crowell group $K_{\s \rho}$ is isomorphic to  $K_\rho$. \end{prop}

\begin{proof} Consider the assignment ${}^iu \mapsto {}^i \mu(u)$, for all $u \in K$, which induces an automorphism of the free group on $\{{}^iu\mid u \in K, i \in [1,N]\}$. The relation ${}^i(uv) = {}^iu\  {}^{i \s\rho(u)}v$ in $K_{\s \rho}$ is mapped to the relation ${}^i(\mu(u)\mu(v)) = {}^i\mu(u)\ {}^{i(\rho(\mu(u))}\!\mu(v)$ in $K_\rho$. Hence there is an induced homomorphism $\mu_\rho: K_{\s\rho} \to K_\rho$. An inverse homomorphism is obtained by replacing $\mu$ with $\mu^{-1}$. Hence $\mu_\rho$ is an isomorphism. \end{proof}

\begin{remark} $\mu_\rho^r$ is an automorphism of $K_\rho$ whenever $\rho$ has period $r$. \end{remark}

\begin{prop}\label{reducible}  If the $K$-orbit of $\N$ consists of $n$ components  with orders $N_1, \ldots, N_n$, then, up to conjugation, $\rho= (\rho_1, \ldots, \rho_n): K \to S_{N_1} \oplus \cdots \oplus S_{N_n}$ and 
$K_\rho \cong K_{\rho_1}* \cdots * K_{\rho_n}$. \end{prop}

\begin{proof} The representation is clearly conjugate to the sum of $\rho_1, 
\ldots, \rho_n$. Consider the associated partition of $\{1, \ldots, N\}$. Any cyclically reduced relator in the presentation (\ref{presentK}) contains superscripts $i$ from just one member of the partition. That $K_\rho$ is the  desired free product follows immediately. 

\end{proof} 

In view of the proposition, we will henceforth assume that the $K$-orbit has only one component. Equivalently, {\sl we assume throughout that $K$ acts transitively on $\N$.} 

When $\rho$ has period $r$, the abelianization $K_\rho^{\rm ab}$ is a finitely generated $\Z[s^{\pm 1}]$-module with $s\  {}^iu= {}^i (\mu^r(u))$, for all $u \in K$. In terms of the presentation (\ref{presentK}), 
$s\ {}^ia_\nu = {}^ia_{\nu+r}, \ldots, s\  {}^ic_\nu = {}^ic_{\nu+r}$.

\begin{definition} The {\it Alexander-Lin} polynomial $D_{\rho, r}(s)$ of a periodic representation $\rho \in \Phi_{S_N}({\cal G})$ with period $r$ is the order $\D_0 K_\rho^{\rm ab}$. \end{definition}

\begin{remark} (1) One should regard $K_\rho^{\rm ab}$ as a twisted Alexander module.  The subscript $r$ is needed in $D_{\rho, r}(s)$ since we do not assume that $r$ is the least period of $\rho$. 

(2) When ${\cal G}$ is the augmented group system of a knot $k$, the invariant $D_{\rho, r}$ coincides with Wada's invariant of the fundamental group of the $r$-fold cyclic cover of $k$.  \end{remark}

The Alexander-Lin polynomial $D_{\rho, r}(s)$ is easily computed using presentation (\ref{presentK}). Form a presentation matrix $J_{\rho,r}$ for $K_\rho^{\rm ab}$ with columns corresponding to ${}^ia_\nu, \ldots, {}^ic_\nu$ and rows corresponding to ${}^ir_\nu, \ldots, {}^is_\nu$, where $ i \in \N$ and $0 \le \nu < r$. Then $D_{\rho, r}(s)$ is the greatest common divisor of the maximal minors of $J_{\rho,r}$. Examples are found below. 

\begin{lemma} \label{shift} If $\rho \in \Phi_{S_N}({\cal G})$ has period $r$, then
$D_{\s\rho, r}(s) = D_{\rho, r}(s)$. \end{lemma} 

\begin{proof} The isomorphism $\mu_\rho: K_{\s \rho} \to K_\rho$ of Proposition \ref{iso} induces a module isomorphism from  
$K_{\s \rho}^{\rm ab}$ to $K_\rho^{\rm ab}$. \end{proof}

When $\rho: K \to S_1$ is the trivial representation, $\D_{\rho, 1}(s)$ is equal to the classical Alexander polynomial $\D(s)$, an invariant of $(G, \e)$. 

\begin{definition} Let $f(t) = c_0 \prod(t-\xi_j) \in {\mathbb C}[t]$. For any positive integer $r$, we set
$f^{(r)}(t) = c_0^r \prod (t-\xi_j^r).$ \end{definition}

\begin{lemma} \label{module} Assume that ${\cal M}$ is a finitely generated torsion $R[t^{\pm 1}]$-module, where $R$ is a Noetherian unique factorization domain. For any positive integer $r$, let ${\cal M}^{(r)}$ denote the same abelian group regarded as an $R[s^{\pm1}]$-module, where $s= t^r.$ If $f(t)$ is the order of ${\cal M}$, then $f^{(r)}(s)$ is the order of ${\cal M}^{(r)}$. \end{lemma}

\begin{proof} Lemma \ref{module} was shown in  \cite{swIJM} in the case that 
$R=\Z$ and ${\cal M}$ has a square presentation matrix. After replacing the field $\C$ with the splitting field of $x^r-1$, regarded as a polynomial over the field of fractions of $R$, the argument of \cite{swIJM} holds for the ring $R$. 

We consider the general case. The module ${\cal M}$ has a square presentation matrix if and only if any finite submodule is trivial (see p. 132 of \cite{hillman}). As in \cite{hillman}, let $z{\cal M}$ be the maximal finite submodule (the {\it pseudonull submodule} of ${\cal M}$). Consider the short exact sequence 
\begin{equation} \label{ses} 0 \to z{\cal M} \hookrightarrow {\cal M} \to {\cal M}/z{\cal M} \to 0. \end{equation}
Since ${\cal M}/z{\cal M}$ has no nontrivial finite submodule, the conclusion of Lemma \ref{module} holds for it. The conclusion also holds for $z{\cal M}$, since its order is necessarily trivial. Regarding (\ref{ses}) as a short exact sequence of $R[s^{\pm 1}]$-modules, the conclusion holds for ${\cal M}$, since its order is the product of the orders of $z{\cal M}$ and ${\cal M}/z{\cal M}$ (see Theorem 3.12 (3) of \cite{hillman}, for example). 
\end{proof} 

\begin{cor} For any positive integer $k$, $D_{\rho, kr}(s) = D_{\rho, r}^{(k)}$. \end{cor}

\begin{prop} \label{divides} If $\rho \in \Phi_{S_N}({\cal G})$ has period $r$, then $\D^{(r)}(s)$ divides $D_{\rho, r}(s)$, where $\D(t)$ denotes the order of $K^{\rm ab}$. 
\end{prop}

\begin{proof} Regard $K^{\rm ab}$ as a $\Z[s^{\pm 1}]$-module with $s=t^r$. By Lemma \ref{module}, its order is $\D^{(r)}(s)$.  Since the projection $p$ of Proposition \ref{quotient} has the property $p \circ \mu_\rho^r = \mu^r \circ p$, it induces a module homomorphism from $K_\rho^{\rm ab}$ onto $K^{\rm ab}$. Hence  $\D^{(r)}(s)$ divides  $D_{\rho, r}(s)$. \end{proof}

\begin{example}\label{bs} Consider the Baumslag-Solitar group $$G=\<x, a \mid x a x^{-1}=a^2\>$$ and epimorphism $\e: G \to \<t \mid \>$ mapping $x \mapsto t, a \mapsto 1$. The kernel $K$ has presentation $\<a_\nu \mid a_{\nu+1}=a_\nu^2\>$, and $\D(t)= t-2$.

There is a period-2 representation $\rho: K \to S_3$ such that $a_{2\nu} \mapsto (123)$ and $a_{2\nu +1} \mapsto (132),$ for all $\nu$.
We regard $S_3$ as a subgroup of ${\rm GL}_3 \Z$ and denote the image $\rho(a_\nu)$ by $A_\nu$. 

Let $s= t^2$. The $\Z[s^{\pm 1}]$-module $K_\rho^{ab}$ has presentation
$$\< {}^i a_0, {}^i a_1 \mid {}^i a_1 = {}^i a_0+ {}^{iA_0}a_0, \ s\ {}^i a_0 = {}^i a_1+ {}^{iA_1}a_1\>.$$ From this we see that 
$$J_{\rho,2} = \begin{pmatrix} -I-A_0 & I \\ s I & -I - A_1\end{pmatrix}.$$

A calculation shows that $D_{\rho, 2}(s) = (s-1)^2(s-4) = (s-1)^2 \D^{(2)}(s)$.  \end{example}

The following will explain the factor $(s-1)^2$ that arises in Example \ref{bs}. 
Consider an augmented group system ${\cal G}=(G, \e, x)$ and a presentation for $G$ of the form (\ref{presentG}). Let $X$ be a standard 2-complex with $\pi_1 X \cong G$ having a single 0-cell $*$, 1-cells corresponding to generators $x, a, b, \ldots, c$, and 2-cells corresponding to relators $r, \ldots, s$. For convenience we label the 1- and 2-cells with the names of corresponding generators and relators. 

Let $X'$ be the infinite cyclic cover of $X$ with $\pi_1 X' \cong K$. It has 0-cells $*_\nu$, 1-cells $x_\nu, a_\nu, \ldots, c_\nu$ and 2-cells $r_\nu, \ldots, s_\nu$ corresponding to the presentation (\ref{presentK}) for $K$. 
{\sl We contract the forest formed by the 1-cells $x_\nu$, but continue to
call the resulting complex $X'$.}

Finally, let $ \widehat {X'}$ be the $N$-sheeted cover of $X'$  corresponding to $\rho$; that is, the cover with $\pi_1 \widehat {X'} \cong {\rm stab}(1).$ (Since $K$ acts transitively on $\N$, by assumption, the various stablizers  ${\rm stab}(i)$ are isomorphic.) For any positive integer $d$, we let $F_d$ denote 
the free group of rank $d$. 

\begin{prop} \label{freeprod} For any $\rho \in \Phi_{S_N}({\cal G})$, the Crowell group $K_\rho$ is isomorphic to the free product $\pi_1 \widehat {X'}*F_{N-1}$. The factor
$\pi_1 \widehat {X'}$ is invariant under the automorphism $\mu_\rho^r$ of $K_\rho$.   \end{prop}

\begin{proof} $K_\rho$ is isomorphic to the fundamental group of $ \widehat {X'}\cup C$, where $C$ is the cone on the 0-skeleton of $\widehat {X'}$.  Transitivity of the action of $K$ on $\N$ enables us to find a cycle $\gamma$ composed of $N$ 1-cells that connect the lifts ${}^i*$. The fundamental group of $C \cup \gamma$ is isomorphic to $F_N$, and $\gamma$ can be chosen to be a member of a basis. By the Seifert-Van Kampen theorem, $K_\rho \cong \pi_1\widehat {X'} *_\gamma F_N \cong  \pi_1\widehat {X'} * F_{N-1}$.  The first statement of the proposition is proved. The second statement is clear. \end{proof}

\begin{cor} \label{homology} If $\rho \in \Phi_{S_N}({\cal G})$ has period $r$, then
$$D_{\rho, r}(s) = \D_0(H_1( \widehat {X'}; \Z) )(s-1)^{N-1}.$$
\end{cor}

\begin{proof}  Let $A$ be a presentation matrix for the $\Z[s^{\pm 1}]$-module 
$H_1( \widehat {X'}; \Z)$ with columns corresponding to generators, and rows to relators.  We assume without loss of generality that the number of rows is at least as great as the the number of columns. The greatest common divisor of the maximal minors of $A$ is $ \D_0(H_1( \widehat {X'}; \Z) )$. 

We can extend $A$ to a presentation matrix for the $\Z[s^{\pm 1}]$-module $K_\rho^{\rm ab}$ with the form 
$$\begin{pmatrix} A & 0 \\ * & D \end{pmatrix},$$  where $D$ is the $(N-1)\times (N-1)$ diagonal matrix ${\rm diag}(s-1, \ldots, s-1)$. 
The greatest common divisor of its maximal minors is $D_{\rho, r}(s)$. 
The proof is completed by noting that any maximal square submatrix with nonvanishing determinant must contain $D$.  \end{proof}

Some periodic representations $\rho \in \Phi_{S_N}({\cal G})$ extend over $G$. This happens if and only if there exists $X \in S_N$ such that $\s\rho(u) = X \rho(u) X^{-1}$. The representation of Example \ref{bs} does not extend. 

When $\rho$ extends, one can consider the twisted Alexander polynomial $\D_\rho(t)$, which is defined in \cite{kl99} as the order of the module $H_1( \widehat {X'}; \Z) $. The following relation follows immediately from Corollary \ref{homology} and Lemma \ref{module}. 

\begin{cor} \label{classicfactor} Assume that  $\rho \in \Phi_{S_N}({\cal G})$ has period $r$. If $\rho$ extends over $G$, then $D_{\rho, r}(s) = \D_\rho^{(r)}(s) \cdot (s-1)^{N-1}$. \end{cor}

\section{3-Manifold groups and reciprocality} Consider an augmented group system ${\cal G} =(G, \e, x)$ such that $G$ is isomorphic to the fundamental group of a compact oriented 3-manifold $X$ with torus boundary. Assume further that $x$ is a peripheral element. Then $\pi_1(\partial X) \cong \Z\oplus \Z$ is generated by $x$ and an element $\ell\in K$. Given $\rho\in\Phi_{S_N}({\cal G})$, the subgroup $\<\ell  \>$ generated by $\ell$ acts by restriction on $\N$. Let $T$ denote the number of orbits. 

\begin{definition} \label{rec} A polynomial $f(s) \in \Z[s^{\pm 1}]$ is {\it reciprocal} if $f(s^{-1})$ is equal to $f(s)$ up to multiplication by units $\pm s^i $. \end{definition}

\begin{theorem}  \label{factors} Assume that  $\rho \in \Phi_{S_N}({\cal G})$ has period $r$.  If $G$ is the group of a 3-manifold $X$ with torus boundary and  $x$ is peripheral, then
\item{(1)} $D_{\rho, r}(s)$ is a reciprocal polynomial;
\item{(2)} $(s-1)^{N+T-2}$ divides $D_{\rho, r}(s)$. \end{theorem}

\begin{proof} Endow $X$ with a cell structure, its 2-skeleton corresponding to a presentation of the form (\ref{presentG}) for $G$. 
Consider the $r$-fold cyclic cover $X_r$ corresponding to $$G=\pi_1X \ {\buildrel\e \over \longrightarrow}\ \Z \to\ \Z/r\Z,$$  a compact 3-manifold with torus boundary that acquires the cell structure lifted from $X$.  The fundamental group $\pi= \pi_1X_r$ is isomorphic to $$\<y, K\mid yuy^{-1} = \mu^r(u)\ (\forall u \in K)\>.$$  Since the representation $\rho: K \to S_N$ has period $r$, it extends over $\pi$,  mapping $y$ trivially. 

The strategy of the proof is to relate $D_{\rho, r}(s)$ to the $\rho$-twisted Reidemeister torsion of $X_r$. We then make use of reciprocality results proved by Kirk and Livingston \cite{kl99}.

The representation $\rho: \pi \to S_N$ induces a right action of $\pi$ on the vector space $V = \C^N$. Consider the cellular chain complex
$$C_*(X_r; V[s^{\pm 1}]_\rho) = (\C[s^{\pm 1}]\otimes_\C V) \otimes_\pi C_*(\tilde X).$$
Here $\tilde X$ is the universal cover of $X$ (and also of $X_r$) with lifted cell structure. The group ring $\Z[\pi]$ acts on the cellular chain complex $C_*(\tilde X)$ on the left by deck transformations. Also, $\C[s^{\pm 1}]\otimes_\C V$ is a right $\Z[\pi]$-module via
$$(p \otimes v) \cdot g = (s^{\e(g)} p) \otimes (v \rho(g)),\ {\rm for}\ g \in \pi.$$ 

Let $X'$ be the infinite-cyclic cover of $X$ (and of $X_r$) with $\pi_1 X' = K$, and let $\widehat X'$ be the $N$-sheeted cover of $X'$ induced by $\rho: K \to S_N$.
By Shapiro's Lemma (see \cite{brown}, for example), 
$H_*(X_r; V[s^{\pm 1}]) \cong H_*(\widehat X'; \C)$.  The $\Z[s^{\pm 1}]$-module structure on the latter group is induced by $\mu^r$. As in \cite{kl99}, the Reidemeister torsion $\tau$ of the complex $C_*(X_r; V[s^{\pm 1}]_\rho)$ is given by 
$${ {\D_0(H_1(\widehat X'; \C))} \over { \D_0(H_0(\widehat X'; \C))\cdot \D_0(H_2(\widehat X'; \C))}   },$$
provided these orders are nonzero.

By Poincar\'e duality, $H_2(\widehat X_r; \C) \cong H_0(\widehat X_r, \partial \widehat X'; \C)$. The latter group is trivial and so $\D_0(H_2(\widehat X'; \C)) = 1$. Also, using the fact that $K$ acts transitively on $\N$, one checks that
$H_0(\widehat X'; \C) \cong \C$, and so $\D_0(H_0(\widehat X'; \C))= s-1$
(the reader can find the straightforward calculations in \cite{swDTAP}).

For the purpose of calculating $H_i(\widehat X'; \C)\cong H_i(\widehat X'; \Z)\otimes_\Z \C$ for $i < 2$, we can ignore 3-cells. If $\D_0(H_1(\widehat X'; \C))$ is the zero polynomial, then by Corollary \ref{homology}, $D_{\rho, r}(s)$ is also zero and the conclusion of the theorem is true trivially. We therefore assume that $\D_0(H_1(\widehat X'; \C))$ is not identically zero, and hence $\tau$ is defined. 

By Theorem 5.1 (and Example 1 of Section 3.3) of  \cite{kl99}, the torsion 
$\tau$ is self-conjugate; that is,
$\bar \tau$ is equal to $\tau$, where $\ \bar{}\ : \sum c_i s^i \mapsto \sum \bar c_i s^{-i}$. Since the denominator $s-1$ has this property, so does
the numerator $\D_0(H_1(\widehat X'; \C))$.  Corollary \ref{homology} then implies that 
$D_{\rho, r}(s) = \bar D_{\rho, r}(s)$, where the polynomials are regarded in $\C[s^{\pm 1}]$.  
However,  $D_{\rho, r}(s)$ has integer coefficients. Hence it is reciprocal.

(2)  Consider the maps induced by the long exact sequence in homology of the pair $(\widehat X', \partial \widehat X')$:
$$H_2(\widehat X',  \partial \widehat X'; \C) \buildrel {\iota} \over {\longrightarrow} H_1( \partial \widehat X'; \C)  \to H_1( \widehat X'; \C).$$
By exactness, $H_1(\partial \widehat X';\C) / {\rm im}\ \iota$ embeds in $H_1(\widehat X'; \C)$. 
The boundary $\partial \widehat X'$ has $T$ components, each homeomorphic to $\S^1 \times {\mathbb R}$, and hence $H_1(\partial \widehat X'; \C) \cong \C^T$. 
Poincar\'e duality implies that $H_2(\widehat X',  \partial \widehat X'; \C)\cong H_0(\widehat X'; \C) \cong \C$. A relative 2-cycle that generates can be represented by a compact orientable surface with boundary comprising $T$ circles onto which $ \partial \widehat X'$ retracts. 
Hence $H_1(\partial \widehat X';\C)/ {\rm im}\ \iota \cong \C^{T-1}$. 
Since the order ideal of a submodule divides that of the module (see, for example, \cite{milnorcover}), $(s-1)^{T-1}$ divides $\D_0(H_1( \widehat {X'}; \C) )$. 
%Arguing as in the last paragraph of the proof of part (1), $(s-1)^{T-1}$ is a factor of $D_{\rho, r}(s)$.  
Hence it also divides $\D_0(H_1( \widehat {X'}; \Z) )$. 
Corollary \ref{homology} ensures  that $D_{\rho, r}(s)$ contains an additional factor of $(s-1)^{N-1}$. 
\end{proof}

\begin{example} \label{7_3} The group $G$ of the 2-bridge knot $k=7_3$ has presentation with generators $x, a$ and a single relator:

$$a^2\cdot x a^{-1}x^{-1}\cdot x^2 a^2 x^{-2}\cdot x^3a^{-1}x^{-3}\cdot x^4 a^2 x^{-4}\cdot x^3 a^{-2} x^{-3} \cdot x^2 a x^{-2}\cdot x a^{-2} x^{-1}.$$
We consider the augmented group system ${\cal G}=(G, \e, x)$, where $\e$ is the abelianization homomorphism mapping $x \mapsto t$.  The commutator subgroup $K$ of $G$ has presentation
$$K= \<a_\nu \mid a^2_\nu a^{-1}_{\nu+1}a^2_{\nu+2}a^{-1}_{\nu+3}a^2_{\nu+4}a^{-2}_{\nu+3}a_{\nu+2} a^{-2}_{\nu+1}\ (\nu\in \Z)\>.$$
From the presentation of $K$, one can read the (untwisted) Alexander polynomial.
$$\D(t) = 2 - 3t+3t^2-3t^3+2t^4.$$
Since the resultant ${\rm Res}(\D(t), t^{13}-1) = 5^4$ vanishes modulo 5, there exist nontrivial elements of $\Phi_{S_{5}}({\cal G})$ with period 13. Such representations are elements of the null space of the $13\times 13$ circulant matrix 
$$\begin{pmatrix} 2 & -3 & 3 & -3 & 2 & 0 & 0 &0   \cdots &0  \\
 0 & 2 & -3 & 3 & -3 & 2 & 0 &\cdots & 0  \\ 
 & & & & \vdots & & &  & \\
 & & & & \vdots & &  & &\\
 -3 & 3 & -3 & 2 & 0 & 0& 0 & \cdots & 2 \end{pmatrix}.$$
 One such representation $\rho$ maps $a_\nu$ ($\nu= 0, \ldots, 12$) respectively to $$\a^4, \a^2, \a, \a, \a^4, \a^4, \a^3, \a, \a^4, 1, 1, 1,\a,$$ where $\a= (12345)$. This representation does not extend to $G$, since $\a$ and $1$ are not conjugate in $S_5$. 
 Calculation shows that $D_{\rho, 13}(s)$ is equal to $(s-1)^8 f(s)g(s)$, where
 $$f(s)= 8192 s^4 - 393 s^3 - 14973 s^2 - 292s + 8192,$$
 and 
 $$g(s)=(64 s^4+ 224 s^3- 801s^2+ 224 s + 64)^2.$$
 
 The longitude $\ell$ is in the second commutator subgroup $[K, K]$, and since the image of $\rho$ is abelian (in fact, cyclic), $\<\ell\>$ acts trivially on $\N$. Theorem \ref{factors} predicts that at least $5+5-2 =8$ factors of $s-1$ will divide $D_{\rho, 13}$. 
 
 The factor $f(s)$ is $\D^{(13)}(s)$, which we expect to see by Corollary \ref{classicfactor}. Each root of $\D(t)$, and hence of this factor, has modulus 1. 
 
 The factor $64 s^4+ 224 s^3- 801s^2+ 224 s + 64$ is not predicted by anything above. Its zeros are all real, equal (approximately) to $1.92246, -5.7693$ and their reciprocals.

\end{example}

\section{Finitely generated kernel $K$} 

\begin{prop} For any $\rho \in \Phi_{S_N}({\cal G})$, the Crowell group $K_\rho$ is finitely generated if and only if $K$ is finitely generated. \end{prop}

\begin{proof} If $K_\rho$ is finitely generated, then so is $K$ by Proposition \ref{quotient}. 

Conversely, if $K$ is generated by finitely many elements $a, b, \ldots, c$, then $K_\rho$ is generated by ${}^ia, {}^ib, \ldots, {}^ic$, where $i$ ranges over $\N$. \end{proof}

\begin{cor} Assume that $\rho \in \Phi_{S_N}({\cal G})$ has period $r$. If $K$ is finitely generated and free, then so is $K_\rho$. In this case, $D_{\rho, r}$ is a monic polynomial of degree $nN$. \end{cor}

\begin{proof} Let $K$ be a free group of rank $n$. Then $K_\rho$ is free of rank $nN$. In this case, $D_{\rho, r}$ is simply the characteristic polynomial of the automorphism of $K_\rho^{\rm ab}$ that is induced by $\mu_\rho^r$, a monic polynomial of degree $nN$.  \end{proof}

\begin{remark} By \cite{koch}, any finitely generated kernel $K$ is free, provided that $G$ has a presentation of deficiency 1. \end{remark}

If $h: H \to H$ is an endomorphism of a finitely generated group $H$, a growth rate ${\rm GR}(h)$ is defined \cite{bowen}. Let $g_1, \ldots, g_n$ be generators for $H$. Let $|| \cdot ||$ denote the word metric on $H$; that is, $||g||$ is the minimum number of occurrences of $g_1^{\pm 1}, \ldots, g_n^{\pm 1}$ in any word describing $g$. Then 
$${\rm GR}(h) = \max_i \limsup_k ||h^k(g_i)||^{1/k}.$$
$GR(h)$ is finite and independent of the generator set. (See \cite{bowen} for this and other basic properties of growth rates.) 

\begin{theorem}\label{bound} Assume that $\rho \in \Phi_{S_N}({\cal G})$ has period $r$. If $K$ is finitely generated, then 
$$(\max\{ |z| \mid D_{\rho, r}(z)=0\})^{1/r} \le {\rm GR}(\mu).$$
\end{theorem}

\begin{proof} 

By Corollary \ref{homology}, it suffices to prove 
$$\max\{ |z| \mid \D_0(H_1(\widehat X'; \Z))(z)=0\}^{1/r} \le {\rm GR}(\mu).$$

Recall that $\mu_\rho^r$ is the automorphism of $K_\rho$ induced by ${}^iu_\nu \mapsto {}^iu_{\nu +r}$, for all $u\in K$. Also, $\pi_1\widehat X'$ is a subgroup of $K_\rho$ invariant under $\mu_\rho^r$, by  Proposition \ref{freeprod}. The restriction $\mu_\rho^r \big\vert_{\pi_1 \widehat X'}$ induces an automorphism of $H_1(\widehat X'; {\mathbb Q})$, and its characteristic polynomial  is $\D_0(H_1(\widehat X'; \Z))$. 

The fundamental group $\pi_1\widehat X'$ is a finite-index subgroup of $K \cong \pi_1 X'$, and the automorphism $\mu^r$ of $K$ restricts to $\mu_\rho^r \big\vert_{\pi_1 \widehat X'}$. 
By a basic property of growth rates, 
${\rm GR}(\mu^r) ={\rm GR}(\mu^r \big\vert_{\pi_1 \widehat X'}).$ 
Also, $|z| \le {\rm GR}(\mu^r)$, for all zeros $z$ of $\D_0(H_1(\widehat X'; \Z))$. The conclusion of the proposition follows from the fact that ${\rm GR}(\mu^r) = {\rm GR}(\mu)^r$.
\end{proof}

\begin{remark} Consider the augmented group system $(G, \e, x)$ comprising the group $G$ of the figure-eight knot $k$, abelianization homomorphism $\e$ and meridianal element $x$. It is well known that the pseudo-Anosov monodromy of $k$ has a stable oriented foliation. It follows that $(3+\sqrt 5)/2$, the larger of the two (real) roots of the the Alexander polynomial $\D(t)=t^2-3t+1$, is equal to the growth rate of $\mu$  (see \cite{ryk}). Hence the bound in Theorem \ref{bound} is the best possible. 
\end{remark}

\begin{cor} \label{twistedbound} Assume that ${\cal G} = (G, \e, x)$ is an augmented group system such that the kernel $K$ of $\e$ is finitely generated. Then for any representation $\rho: G \to S_N$, 

$$\max\{ |z| \mid \D_\rho(z)=0\}\le {\rm GR}(\mu).$$
\end{cor}

\begin{proof} Let $r$ be the order of $\rho(x)$. The representation $\rho$ restricted to $K$ is periodic with order $r$, and the result follows immediately from Corollary \ref{classicfactor}.  \end{proof}

\begin{remark} Theorem 1.3 of \cite{mcmullen} implies that the bound of corollaries \ref{bound} and \ref{twistedbound} is strict whenever $G$ is the group of a hyperbolic 3-manifold  such that the pseudo-Anosov monodromy associated to $\e$ does not have a stable oriented foliation. 
\end{remark}

Let $k$ be a nonfibered knot with augmented group system ${\cal G}$. Consider only representations $\rho \in \Phi_{S_N}({\cal G})$ of period $r$ for which $D_{\rho, r}(s)$ is nonzero.  In view of Theorem \ref{bound} we ask: 

\begin{question} Does $\max \{ |z| \mid D_{\rho, r}(z) = 0\}^{1/r}$ have a bound independent of both $\rho$ and $r$? \end{question}

\section{The curious case of the vanishing polynomial} 

When $K$ is not finitely generated, it can happen that $D_{\rho, r}(s)$ is identically zero. We explain how this can happen. 

Given an augmented group system $(G,\e, x)$, the group $G$ can be described as an HNN decomposition 
\begin{equation} \label{HNN} G = \<B, x \mid x a x^{-1} = \phi(a)\ (\forall a \in U)\>,\end{equation}
where $B$ is a finitely generated subgroup of $K$,  $U$ and $V$ 
are finitely generated subgroups of $B$, and $\phi: U \to V$ is an isomorphism (see \cite{ls}, for example). 
Let $B_\nu = x^\nu B x^{-\nu}$, $U_\nu = x^\nu U x^{-\nu}$ and $V_\nu = x^\nu V x^{-\nu}$, for all $\nu \in \Z$. Note that $\phi$ induces an isomorphism between $V_\nu$ and $U_{\nu+1}$. Then $K$ is an infinite free product with amalgamation, 
$$K = \cdots *_U B *_U B *_U\cdots,$$ 
where $\phi$ identifies $V_\nu \subset B_\nu$ with $U_{\nu+1}\subset B_{\nu+1}$.

A basic property of amalgamated free products implies that $B_\nu$ is a subgroup of $K$, for each $\nu$. 
The subgroups $U_\nu$ of $K$ act on $\N$ by restriction. 

\begin{prop}\label{vanish} Assume that $\rho \in \Phi_{S_N}({\cal G})$ has period $r$. If $U_\nu$ does not act transitively on $\N$ for some $\nu\in \{ 1, \ldots, r-1\}$, then $D_{\rho, r}(s) = 0$. \end{prop}

\begin{proof} Our argument was motivated by ideas of \cite{fv08}. 

Assume that $U_\nu$ does not act transitively on $\N$, for some $\nu$.  By Lemma \ref{shift}, we can assume that $\nu=0$. 

As in the proof of Theorem \ref{factors}, consider the finite-index subgroup $\< y, K \mid y u y^{-1} = \mu^r(u) \ (u \in K)\>$ of $G$, with $y= x^r$. 
It has an HNN decomposition
$$\<B_{[1,r]}, y \mid y b y^{-1} = \mu^r(b)\ (\forall b \in U)\>,$$  with base $B_{[1,r]}$ equal to the subgroup of $K$ generated by $B_1, \ldots, B_r$.

Let $Y$ be a cell complex such that $\pi_1 Y \cong B_{[1,r]}$, with subcomplexes $Y', Y''$ such that $\pi_1Y' \cong U$ and $\pi_1 Y'' \cong V$. A complex  with fundamental group $K$ is formed from countably many copies $(Y_\nu; Y'_\nu, Y''_\nu)$ by gluing end to end according to the amalagamating map $\phi$. 

The representation $\rho$ induces an $N$-sheeted cover $\widehat Y$. It restricts to each $Y_\nu, Y'_\nu, Y''_\nu$, inducing respective $N$-sheeted covers $\widehat Y_\nu; \widehat Y'_\nu, \widehat Y''_\nu$ such that 
$\widehat Y$ is the union of the $\widehat Y_\nu$ glued end to end. 

The Mayer-Vietoris long exact sequence contains
$$ \cdots \to K_\rho^{\rm ab} \to \oplus_\nu H_0(\widehat Y'; \Z) \to \oplus_\nu H_0(\widehat Y; \Z) \to \cdots$$ Since $K$ acts transitively on $\N$, so does $B_{[1,r]}$. Hence each $H_0(\widehat Y_\nu; \Z) \cong \Z$, and  
$\oplus_\nu H_0(\widehat Y; \Z) \cong \Z[s^{\pm 1}]$. However, the $U$-orbit of $\N$ consists of $q$ components with $q>1$. Hence each $H_0(\widehat Y'_\nu; \Z) \cong \Z^q$, and $\oplus_\nu H_0(\widehat Y'; \Z) \cong (\Z[s^{\pm 1}])^q$. It follows from the exact sequence that $K_\rho^{\rm ab}$ cannot be a $\Z[s^{\pm 1}]$-torsion module. Thus $D_{\rho, r}(s)$ is zero.

\end{proof}

When $G$ is the group of a knot $k$, the fundamental group of $\S^3$ split along any incompressible Seifert surface $S$ for $k$ is an HNN base $B$  for $G$. The amalgamating subgroups are $U=\pi_1 S$ and $V=x (\pi_1 S )x^{-1}$.
Proposition \ref{vanish} implies that if $\pi_1 S$ does not act transitively on $\N$, then $D_{\rho, r}(s)$ must be identically zero. However, $\pi_1 S$ can fail to act transitively only if $k$ is nonfibered, since otherwise $S^3$ split along $S$ is a product, causing $B, U$ and $V$ to coincide. 

In \cite{swBLM}  we conjectured that if $k$ is nonfibered, then there exists a representation $\rho: K \to S_N$, for some $N$, such that $\rho(U)$ is a proper subgroup of $\rho(B)$. In such a case, it is not hard to find another representation, possibly with larger $N$, such that $B$ acts transitively on $\N$ but $U$ does not (see \cite{swJPAA}).  Our conjecture then is equivalent to the following.

 \begin{conj} \label{polyzero} Let $k$ be a knot with augmented group system ${\cal G} = (G, \e, x)$. Then $k$ is nonfibered if and only if for some integer $N>0$ and periodic representation $\rho \in \Phi_{S_N}({\cal G})$, the polynomial $D_{\rho, r}(s)$ is zero.
 \end{conj}
 
A group $G$ is {\it subgroup separable} if, given any finitely generated subgroup $U$ and element $a \notin U$, there exists a homomorphism $\g$ from $\pi$ to a finite group $\Sigma$ such that $\g(a) \notin \g(U)$. 
In \cite{wi} D. Wise announced work, some of it joint with F. Haglund, C. Hruska, T. Hsu and M. Sageev, showing that a hyperbolic 3-manifold with an incompressible geometrically finite surface has a subgroup separable fundamental group.  Combined with earlier results of F. Bonahon \cite{bon} (or alternatively \cite{agol}, \cite{cg}) showing that  an incompressible Seifert surface of a hyperbolic knot $k$ is either a virtual fiber or else is geometrically finite, the  results of \cite{wi} would imply that Conjecture \ref{polyzero} holds for all hyperbolic knots. 

%The following result is included for the sake of completeness. It follows from subgroup separability and Theorem 3.4 of \cite{swBLM}
%(see also \cite{swNFRS}).  Theorem 1 of Friedl and Vidussi \cite{fv08} is a general statement for irreducible 3-manifolds with the appropriate subgroup separability. 

%
%\begin{theorem} \label{hyper} A hyperbolic knot is fibered if and only if $\D_\g$ is monic for every finite-image representation. \end{theorem}

%Theorem \ref{hyper} can be viewed as a generalization of K. Murasugi's theorem that any alternating knot is fibered if and only if its classical Alexander polynomial is monic. 
 
\section{Bounding the degree of $D_{\rho, r}$}

An HNN decomposition 
of the form (\ref{HNN}) provides a lower bound for the degree of $D_{\rho, r}$. 

\begin{theorem} \label{Bound} Let ${\cal G} = (G, \e, x)$ be an augmented group system and assume that $\rho: K \to S_N$ is a representation of period $r$. Let $$G = \<B, x \mid x a x^{-1} = \phi(a)\ (\forall a \in U)\>$$ be an HNN decomposition for $G$. Then 
$$\deg D_{\rho, r} \le N\cdot {\rm rk}(U),$$
where ${\rm rk}(U)$ is the minimum number of elements needed to generate $U$. \end{theorem}

\begin{proof} The theorem is a special case of Proposition 2.4 of \cite{preprint}. The $\Z[s^{\pm 1}]$-module $K_\rho^{\rm ab}$ is isomorphic to an infinite free abelian product 
$$ \cdots \oplus_{\hat U} \hat B \oplus_{\hat U} \hat B \oplus_{\hat U}\cdots$$
with identical amalgamations given by monomorphisms $\hat U\ {\buildrel g  \over \leftarrow}\  \hat B\  {\buildrel f \over \rightarrow}\  \hat U$, where $\hat U$ and $\hat B$ are finitely generated abelian groups and $\hat U$ is generated by no more than ${\rm rk}(U)$ elements. (The monomorphisms are necessarily isomorphisms when $K_\rho$ and hence $K_\rho^{\rm ab}$ are finitely generated.) As in \cite{preprint}, one constructs a presentation matrix for $K_\rho^{\rm ab}$ with $Q$ columns corresponding to generators of $\hat B$ and $P+q$ rows corresponding to the $P$ defining relations for $\hat B$ as well as the $q$ amalgamating relations of the form $ s g(u) = f(u)$. Here $u$ ranges over $q$ generators for $\hat U$. Without loss of generality, we can assume that $Q \ge q$ and $P+q \ge Q$. Then
$D_{\rho, r}$ is the greatest common divisor of the $Q \times Q$ minors. The conclusion of the theorem follows immediately. 
\end{proof}

\begin{cor} \label{corbound} If $G$ is the group of a knot $k$, then for any representation $\rho: K \to S_N$ of period $r$, 
$$\deg D_{\rho, r} \le N \cdot 2\ {\rm genus}(k),$$ 
with equality for fibered knots. 
\end{cor}

\begin{remark}  If $G$ is the group of a knot $k$, the representation $\rho$ is trivial, and $U$ is the fundamental group of a Seifert surface of minimal genus $g$, then the statement of Theorem \ref{Bound}
reduces to the well-known result that the degree of the untwisted Alexander polynomial of $k$ is no greater than $2 g$. \end{remark} 

\begin{question}  Does the group of some knot admit an HNN decomposition (\ref{HNN}) such that the rank of $U$ is less than twice the genus of $k$? \end{question}

\section{Mahler measure of $D_{\rho, r}$}  \begin{definition} The Mahler measure of a nonzero polynomial $p(s) = c_d s^d + \cdots + c_1 s + c_0 \in {\mathbb C}[s]$, with $c_d \ne 0$,  is 
$$M(p)=|c_d| \prod_{i=1}^d \max\{|\lambda_i|, 1\},$$ where $\lambda_1, \ldots, \lambda_d$ are the roots of $p(s)$. \end{definition}

Mahler measure is clearly multiplicative, and the Mahler measure of a unit is 1. Hence we can extend the definition of Mahler measure to Laurent polynomials that are well defined up to multiplication by units. In particular,
the Mahler measures of nonzero Alexander-Lin polynomials are well defined. 

It is well known that for polynomials $p$ with integer coefficients, $M(p)=1$ if and only if $p$ is equal up to a unit factor to a product of cyclotomic polynomials. (See \cite{ew} for this as well as other basic facts about Mahler measure.) 

In \cite{swDTAP} the authors characterized the  Mahler measure of twisted Alexander polynomials for knots in terms of the growth of homology torsion in cyclic covers. The results apply to a large class of knot group representations that include finite permutation representations, and they generalize previous results \cite{swTOP2} for untwisted Alexander polynomials. 

When periodic representations of the commutator subgroup are considered, a characterization of the Mahler measure of the resulting twisted Alexander-Lin polynomial is simpler since, as in \cite{swTOP2}, it can stated in terms of {\sl branched} cyclic covers. 

Assume that ${\cal G} = (G, \e, x)$  is an augmented group system associated to a knot $k$, and  $\rho \in \Phi_{S_N}({\cal G})$ is a representation with period $r$. 
Then $\rho$ induces a representation of the $rn$-fold cyclic cover $M_{rn}$ of 
$\S^3$ branched over $k$, for any positive integer $n$. Let $\widehat M_{rn}$ denote the $N$-fold cover induced by $\rho$. 
We define the {\it twisted torsion number} $b_{\rho, rn}$ to be the order of the torsion subgroup of $H_1 (\widehat M_{rn}; \Z)$.

\begin{theorem}\label{mahler} If $D_{\rho, r}(s) \ne 0$, then 
$$M(D_{\rho, r}(s))= \lim_{n\to \infty} (b_{\rho, rn})^{1/n} .$$ \end{theorem}

\begin{proof} Construct a standard 2-complex $X$ with fundamental group $G$, as in the proof of Corollary \ref{homology}. As in the proof, let $X'$ denote the infinite cyclic cover with $1$-cells $x_\nu$ contracted to a single $0$-cell $*$, and let $\widehat X'$ be the $N$-fold covering space of $X'$ corresponding to $\rho$. One sees using the Wang homology sequence that $H_1(\widehat M_{rn}; \Z)$ is isomorphic to $H_1(\widehat X'; \Z)/(s^n-1)$.  By the corollary, it suffices to prove the result with 
$D_{\rho, r}$ replaced by the order $\D_0(H_1(\widehat X'; \Z))$.
Since $G$ has a presentation with one more generator than the number of relators, the $\Z[s^{\pm 1}]$-module $H_1(\widehat X'; \Z)$ has a square matrix presentation. Theorem 2.10 of \cite{swTN}, which expresses the order of any $\Z[s^{\pm 1}]$-module with a square presentation matrix as the exponential growth rate of similarly defined torsion numbers, completes the proof. \end{proof}

We conclude with a theorem that illustrates the advantage of considering representations of the augmentation subgroup $K$ rather than those of $G$. 

\begin{theorem} \label{bounds} Assume that $(G, \e, x)$ is an augmented group. Then there exists a positive number $B$ such that for any period-$r$ finite-image abelian representation $\rho$, 
$$\max \{ |z| \mid D_{\rho, r} (z)=0,\  D_{\rho, r} \ne 0\}^{1/r} \le B.$$
\end{theorem}

\begin{proof} The representation $\rho$ decomposes as the sum of 1-dimensional unitary representations $\rho_1, \ldots, \rho_N$. Moreover $D_{\rho,r} = D_{\rho_1,1}  \cdots D_{\rho_N,1} $. Hence it suffices to consider the case $N=1$. 

The group $G$ has a presentation of the form $\<x, a_1, \ldots, a_n \mid r_1, \ldots, r_m\>$, where each $a_i$ is in the kernel $K$. As usual, we assume without loss of generality that $n \le m$. 
Let $M$ be the maximum length of any relator. 

Consider the presentation  $\< a_{i, \nu} \mid r_{j, \nu}\>$ for $K$, as in Section 2,  
where $a_{i, \nu}$ and $r_{j, \nu}$ denote $x^\nu a_i x^{-\nu}$ and $x^\nu r_j x^{-\nu}$, respectively.  In the $m r \times n r$ presentation matrix $J_{\rho, r}$ of $K_\rho^{\rm ab}$  (see Section 4), the maximum number of nonzero entries in any row is $M$, and no entry consists of more than $M$ monomials.  Any maximal minor is a polynomial with coefficient sum no greater than $M^{2nr}$. 
The Mahler measure of such a polynomial is no greater than $M^{2nr}$ (see \cite{ew}), and hence no zero has modulus greater than $M^{2nr}$. Since $D_{\rho, r} $ is the greatest common divisor of maximal minors of $J_{\rho,r}$, the desired conclusion follows with $B= M^{2n}$.
\end{proof} 

\begin{cor} Assume that $(G, \e, x)$ is an augmented group. Then there exists a positive number $B$
such that for any finite-image metabelian representation $\rho$ of $G$,  
$$\max \{ |z| \mid \D_\rho (z)=0,\  \D_\rho \ne 0\} \le B.$$
 \end{cor}

\begin{proof} The representation $\rho$ restricts to a finite-image abelian representation of $K$ with period $r$, where $r$ is the order of $\rho(x)$. By Corollary \ref{classicfactor}, $\D_\rho^{(r)}$ divides $D_{\rho, r}$. The desired result follows immediately from Theorem \ref{bounds}. 

\end{proof}

\begin{question} Can a bound $B$ in Theorem \ref{bounds} be found that is independent of the augmentation $\e$? 
\end{question}

\end{document}